\documentclass[12pt]{amsart}

\usepackage{amsfonts,amssymb,eucal}
\usepackage{amsthm}
 \usepackage{amsmath} 
\usepackage{amscd}
 \usepackage{latexsym} 
\input{cyracc.def}

\numberwithin{equation}{section}
\newcommand{\SL}{\text{SL}_2(\R )}
\newcommand{\Ad}{\ensuremath{{\mbox{\rm{Ad}}}}}
\newcommand{\ad}{{\mbox{\rm ad}}}
\newcommand{\e}{{\mbox{\rm e}}}

\newcommand{\mb}[1]{{\mbox{\boldmath{$#1$}}}}
\newcommand{\mc}[1]{{\mathcal{#1}}}
\newcommand{\got}[1]{{\mathfrak{#1}}}
\newcommand{\db}[1]{{\mathbb{#1}}}
\newcommand{\gata}{\blacksquare}
\newcommand{\pa}{\partial}

\newcommand{\R}{\ensuremath{\mathbb{R}}}
\newcommand{\C}{\ensuremath{\mathbb{C}}}

\newcommand{\Hi}{\ensuremath{\mathcal{H}}}

\newtheorem{Remark}{Remark}

\newtheorem{Proposition}{Proposition}
\newtheorem{lemma}{Lemma}
\newtheorem{corollary}{Corollary}

\theoremstyle{definition}

 \def\i{\mathrm{i}}

\begin{document}
\title [The fundamental conjecture for the Siegel-Jacobi disk]{Consequences of the fundamental conjecture for the motion
  on the  Siegel-Jacobi disk}
\author{Stefan  Berceanu}
\address[Stefan  Berceanu]{Horia Hulubei National
 Institute for Physics and Nuclear Engineering\\
        Department of Theoretical Physics\\
 P.O.B. MG-6, 077125 Magurele, Romania}
\email{ Berceanu@theory.nipne.ro}

\begin{abstract} 
We find the homogenous K\"ahler diffeomorphism  $FC$ which expresses the K\"ahler two-form on the
Siegel-Jacobi domain  $\mc{D}^J_1=\C\times \mc{D}_1$ 
as the sum of the K\"ahler two-form on $\C$ and the one on the Siegel ball
$\mc{D}_1$. The classical motion and quantum evolution  on
$\mc{D}^J_1$ determined by a linear Hamiltonian in the generators of
the Jacobi group $G^J_1=H_1\rtimes\text{SU}(1,1)$ is described by a
 Riccati equation on $\mc{D}_1$ and a linear first order
differential equation in $z\in\C$, where $H_1$ denotes the 
3-dimensional Heisenberg group. When  the transformation $FC$
is applied,  the first order differential  equation for the variable $z\in \C$ decouples of the
motion on the Siegel disk. Similar considerations are presented for
the Siegel-Jacobi space $\mc{X}^J_1=\C\times\mc{X}_1$,
where $\mc{X}_1$ denotes the Siegel upper half plane. 
\end{abstract}
\maketitle

\noindent

\setcounter{tocdepth}{3}


\section{Introduction}
The Jacobi group   \cite{ez}   is the semidirect product 
$G^J_n=H_n\rtimes\text{Sp}(n,\R)$, 
where $H_n$ denotes the real
$(2n+1)$-dimensional Heisenberg group.  Several generalizations are
known \cite{TA99}, 
\cite{LEE03}. The   homogenous   K\"ahler Siegel-Jacobi 
domains $\mc{D}^J_n=\C^n\times\mc{D}_n$, where $\mc{D}_n$ is the
Siegel ball,  are nonsymmetric domains associated to the nonreductive  Jacobi groups
by the generalized Harish-Chandra embedding \cite{satake}, \cite{LEE03}, 
\cite{Y02}-\cite{Y08}. 
The holomorphic irreducible unitary representations of the Jacobi groups
based on Siegel-Jacobi domains have been constructed \cite{bb,bs,TA90,TA99,gem}.

Some coherent state systems based on Siegel-Jacobi domains have been
investigated in the framework of quantum mechanics, geometric quantization,
dequantization, quantum optics, nuclear structure, and signal processing 
\cite{KRSAR82, Q90, SH03,jac1,sbj,sbcg}.  The Jacobi group was investigated by
physicists under other names as 
 {\it Hagen} \cite{hagen},  {\it Schr\"odinger} \cite{ni},    or {\it Weyl-symplectic} group
 \cite{kbw1}. The Jacobi
group is  responsible  for the  {\it squeezed states} \cite{ken,stol,lu,yu,ho} in quantum optics
\cite{mandel,ali,siv,dr}.

In the papers \cite{jac1,sbj}  the  Jacobi group is studied in
connection with the group-theoretic approach to coherent states 
\cite{perG}. We have attached  to the Jacobi group
$G^J_n$ coherent states  based on
Siegel-Jacobi disk  $\mc{D}^J_n$ \cite{sbj}. In this paper we consider the case
of the Jacobi group $G^J_1$  \cite{jac1}. We have determined the
$G^J_n$-invariant Kh\"aler two-form $\omega_n$ on $\mc{D}^J_n$ \cite{sbj} from  the K\"ahler potential, and,  with the partial Cayley
transform, we have determined the K\"ahler two-form $\omega'_n$  on the
Siegel-Jacobi upper-half plane $\mc{X}^J_n=C^n\times\mc{X}_n$, where
$\mc{X}_n$ is the Siegel upper half-plane. The K\"ahler two form
$\omega'_1$  on $\mc{X}^J_1=C\times\mc{X}_1$  was
firstly investigated by K\"ahler himself \cite{cal3} and Berndt \cite{bern}, while  $\omega_n$
and  $\omega'_n$ have been investigated also by Yang \cite{Y07},\cite{Y08}. $\omega_1$ is
the sum of two terms, one $\omega_{\mc{D}_1}$ describing the the K\"ahler  two form on 
$\mc{D}_1$,  the other one is 
$(1-w\bar{w})^{-1}A\wedge\bar{A}$,  where $ A=dz +\bar{\eta}dw$, and 
$\eta=(1-w\bar{w})^{-1}(z+\bar{z}w)$,  $z\in\C, w\in \mc{D}_1$ \cite{jac1}. Let as
denote by $FC$ the change of variables $FC: 
\mc{D}^J_1=\C\times\mc{D}_1\ni (z,w)\rightarrow (\eta,w)$.  It is
find out that if we express $\omega_1$ in the coordinates $(\eta, w)$,
then the K\"ahler
two-form
becomes $    \omega_0: =  FC^*(\omega_1)=\omega_{\mc{D}_1}+\omega_{\C}$,
invariant to the action of $G^J_1$ on the product manifold
$\C\times\mc{D}_1$.  We put
this change of variables in connection with the celebrated fundamental
conjecture  of  Gindikin-Vinberg \cite{GV},\cite{DN} on the
homogeneous K\"ahler Siegel-Jacobi disk. Similar considerations
are presented for the homogenous K\"ahler Siegel-Jacobi space $\mc{X}^J_1$. 

{\it Dequantization} in the
Berezin approach \cite{berezin2},\cite{berezin1}  of a dynamical system problem with Lie group of
symmetry $G$   on a Hilbert space
$\Hi$ in the simple case of linear Hamiltonian was  considered  in 
\cite{sbcag},\cite{sbl}. Linear Hamiltonians in generators of the Jacobi group appear in many
physical problems of quantum mechanics, as in the case of the  quantum oscillator
acted on by a variable external force \cite{fey}, \cite{sw},
\cite{hs}. The same problem was considered in the group-theoretic
approach in \cite{pp,ppAV,perG}, where linear
Hamiltonians in the generators of $\text{SU}(1,1)$ or $\text{Isp}(2)$
are considered. A similar treatment has been used in the case of  quantum
dynamics of trapped ions \cite{viorica}.

The paper is laid  out as follows. \S \ref{JAC11}, devoted to
$G^J_1$, starts  in \S\ref{jac1} with the Jacobi algebra $\got{g}^J_1$
and
Perelomov's coherent states defined on $\mc{D}^J_1$. Several facts
concerning a holomorphic representation of $\got{g}^J_1$ as first order
differential operators with polynomial coefficients on $\mc{D}^J_1$
are summarized in  Lemma \ref{mixt}, which is essential for \S \ref{CLSQ}.
 Remark \ref{rem31} is new.  In\S \ref{REALGR} we recall some facts referring to the
real Jacobi group $G^J_1(\R)$, firstly considered by K\"ahler and
Berndt  \cite{cal1,cal2,cal3,bern,bs}.  Proposition \ref{cevan}   and
Remark \ref{actnea1} are  extracted  from
\cite{jac1}.  In \S
\ref{ok} we find out the homogenous K\"ahler isomorphism $FC:$ 
$\mc{D}^J_1\rightarrow \C\times\mc{D}_1$
($FC_1:~ \mc{X}^J_1\rightarrow\C\times\mc{X}_1$, respectively).  In \S \ref{CLSQ}  we study classical and quantum motion
determined by linear Hamiltonians in generators of the Jacobi group   $G^J_1$  on
 Siegel-Jacobi
disk and upper half plane. We calculate   the Berry 
phase and the dynamical phase \cite{swA} on $\mc{D}^J_1$ 
 determined
by a 
Hamiltonian  linear in the generators of the Jacobi group. The simple
example of an autonomous system is explicitly solved. 
The main results   of the paper are stated  in
Proposition \ref{MBNM},  Corollary \ref{siasta}, and Proposition
\ref{EQLIN1}.

\section{The Jacobi group $G^J_1$}\label{JAC11}

 \subsection{The Jacobi algebra  $\got{g}^J_1$ and Perelomov's
   coherent states }\label{jac1}
We consider a realization of the 3-dimensional    Heisenberg    
Lie algebra 
\begin{equation}\label{nr0}\got{h}_1\equiv
<\i s 1+xa^{\dagger}-\bar{x}a>_{s\in\R ,x\in\C} ,\end{equation}
 where $ a^{\dagger}$ (${ a}$) are  the boson creation
(respectively, annihilation)
operators,   $ [a,a^{\dagger}]=1$.

 The Jacobi algebra is defined  as the  the semi-direct sum
$\got{g}^J_{1}:= \got{h}_1\rtimes \got{su}(1,1)$, where
\begin{equation}\label{nr}
\got{su}(1,1)=
<2\i\theta K_0+yK_+-\bar{y}K_->_{\theta\in\R ,y\in\C} ,
 \end{equation} 
\begin{equation}\label{nr1}
  \left[ K_0, K_{\pm}\right]=\pm K_{\pm}, ~
\left[ K_-,K_+ \right]=2K_0; 
\end{equation} 
\begin{subequations}\label{baza}
\begin{eqnarray}
 \left[a,K_+\right] & = & a^{\dagger},\left[ K_-,a^{\dagger}\right] =a, 
\left[ K_+,a^{\dagger}\right]=\left[ K_-,a \right]= 0,\label{baza1}\\
   \left[ K_0  ,a^{\dagger}\right] &  =
  &\frac{1}{2}a^{\dagger},~ 
\left[ K_0,a\right]
= -\frac{1}{2}a .\label{baza2}
\end{eqnarray}
\end{subequations}



Let us  suppose that we know the  derived representation $d\pi$ of 
the Lie algebra $\got{g}^J_{1}$ 
of the Jacobi group $G^J_{1}$. For a Lie algebra $\got{g}$, if $X\in\got{g}$, we denote $\mb{X}=d\pi (X)$.
 We impose to the cyclic vector $e_0$ to verify simultaneously
 the conditions
\begin{equation}\label{cond}
\mb{a}e_0  =  0, ~
 {\mb{K}}_-e_0  =  0,~
{\mb{K}}_0e_0  =  k e_0;~ k>0, 2k=2,3,... , 
\end{equation}
and we have considered in  the last relation in (\ref{cond}) the positive  discrete series
representations $D^+_k$ of $\text{SU}(1,1)$  \cite{bar47}.

 Perelomov's coherent state   vectors   associated to the group $G^J_{1}$ with 
Lie algebra the Jacobi algebra $\got{g}^J_{1}$, based on Siegel-Jacobi
disk $ \mathfrak{D}^J_{1}  =  H_1/\R\times \text{SU}(1,1)/\text{U}(1)$
$ =\C\times\mc{D}_1$,  
are defined as 
\begin{equation}\label{csu}
e_{z,w}:=\e^{z\mb{a}^{\dagger}+w{\mb{K}}_+}e_0, ~z,w\in\C,~ |w|<1 .
\end{equation}

The formulas below are obtained \cite{jac1} using the relation  $\Ad(\exp X)=\exp(\ad_X)$
\begin{subequations}\label{summa1x1}
\begin{eqnarray}
& & \mb{a}^+e_{z,w} =\frac{\pa}{\pa z} e_{z,w};~\mb{a} =(z+w\frac{\pa}{\pa z}) e_{z,w};\\
 & & \db{K}_+ e_{z,w}=\frac{\pa}{\pa w} e_{z,w};~\db{K}_0=(k+\frac{1}{2}z\frac{\pa}{\pa z}+
w\frac{\pa}{\pa w})e_{z,w};\\
& & \db{K}_-e_{z,w}=(\frac{1}{2}z^2+2kw +zw\frac{\pa}{\pa z}+
w^2\frac{\pa}{\pa w})e_{z,w} .
\end{eqnarray}
\end{subequations}

With (\ref{summa1x1}), the general scheme \cite{sbcag,sbl} associates to elements of the
Lie algebra $\got{g}$,  first order holomorphic differential operators
with polynomial coefficients
 $X\in\got{g}\rightarrow\db{X}$:  
\begin{lemma}\label{mixt}The differential action of the generators
 of
the Jacobi algebra {\em (\ref{baza})} is given by the formulas:
\begin{subequations}\label{summa}
\begin{eqnarray}
& & \mb{a}=\frac{\pa}{\pa z};~\mb{a}^+=z+w\frac{\pa}{\pa z} ,
~z,w\in\C, ~|w|<1; \\
 & & \db{K}_-=\frac{\pa}{\pa w};~\db{K}_0=k+\frac{1}{2}z\frac{\pa}{\pa z}+
w\frac{\pa}{\pa w};\\
& & \db{K}_+=\frac{1}{2}z^2+2kw +zw\frac{\pa}{\pa z}+w^2\frac{\pa}{\pa
w} . 
\end{eqnarray}
\end{subequations}
\end{lemma}

Acting on $e_{z,w}$, the differential operators (\ref{summa1x1}) are
not independent: 
\begin{Remark}\label{rem31}
Perelomov's coherent state vector \emph{(\ref{csu})} verifies the
  system of differential equations $X e_{z,w}=0$, $Y e_{z,w}=0$, where
\begin{eqnarray}\label{XYVG}
X & =  & k\mb{a}+(\frac{z^2}{2}-kw)\mb{a}^++wk\db{K}_+-z\db{K}_0,   \\
Y &=  &\frac{1}{2}z^3\mb{a}^+ + (2kw+z^2)w\db{K}_+-(4kw+z^2)\db{K}_0+2k\db{K}_-.
\end{eqnarray}
\end{Remark}
{\it Proof.}  We want to determine the constants $A-E$ such that
$$(A\mb{a}+B\mb{a}^{\dagger}+C\db{K}_++D\db{K}_0+E\db{K}_-)e_{z,w}=0. $$
We use the relations (\ref{summa1x1}), and equate with 0 the
coefficients of $constant$, $\frac{\pa}{\pa z}$, and $\frac{\pa}{\pa
  w}$. We express four of the constants $A-E$ as function of
other two. Choosing $A$ and $E$ as independent, we get  the expressions
given in (\ref{XYVG}). \hfill 
  $\gata$

We consider the  displacement operator
\begin{equation}\label{deplasare}
D(\alpha ) =\exp (\alpha \mb{a}^{\dagger}-\bar{\alpha}\mb{a})=\exp(-\frac{1}{2}|\alpha
|^2)  \exp (\alpha \mb{a}^{\dagger})\exp(-\bar{\alpha}\mb{a}), 
\end{equation}
and let  us denote by $S$  {\it the unitary squeezed operator}  --
 the $D^k_+$ representation of the group $\text{SU}(1,1)$. 
We introduce the notation $\underline{S}(z)=S(w)$, $z,
w\in\C, |w|<1$,  where
\begin{subequations}
\begin{eqnarray}
\underline{S}(z) & = &\exp (z{\mb{K}}_+-\bar{z}{\mb{K}}_-)=\exp (w{\mb{K}}_+)\exp (\eta
{\mb{K}}_0)\exp(-\bar{w}{\mb{K}}_-) , \label{u2}\\
w & = &  \frac{z}{|z|}\tanh \,(|z|),~~ \label{u3}
\eta= \log (1-w\bar{w}) .\label{u5}
\end{eqnarray}
\end{subequations}
 We introduce also the normalized    ({\it squeezed})     CS vector  
$\Psi_{\alpha, w}:=D(\alpha )S(w) e_0$  \cite{stol} . 

We introduce the auxiliary  operators  \cite{jac1}:
\begin{equation}\label{ssq}
{\mb{K}}_+  = 
 \frac{1}{2}(\mb{a}^{\dagger})^2+{\mb{K}}'_+ ,~
{\mb{K}}_-  =  \frac{1}{2}(\mb{a}^{\dagger})^2+{\mb{K}}'_- ,~
{\mb{K}}_0  = \frac{1}{2}(\mb{a}^{\dagger}\mb{a}+\frac{1}{2})+{\mb{K}}'_0 ,
\end{equation}
which have the properties
\begin{gather}
{\mb{K}}'_-e_0  =  0 ,
{\mb{K}}'_0e_0 =  k'e_0; ~k=k'+\frac{1}{4};\\
\label{kkk3} 
 [{\mb{K}}'_{\sigma},\mb{a}]=[
{\mb{K}}'_{\sigma},\mb{a}^{\dagger}]=0,~\sigma =\pm ,0 , 
 [ {\mb{K}}'_0,{\mb{K}}'_{\pm}]=\pm {\mb{K}}'_{\pm};
[ {\mb{K}}'_-,{\mb{K}}'_+]=2{\mb{K}}'_0 .
\end{gather}

We recall the orthonormal  system of coherent states associated to the
group
$\text{SU}(1,1)$: 
\begin{equation}\label{ort1}
e_{k,k+m}:=a_{km}({\mb{K}}_+)^me_{k,k};~~a^2_{km}=\frac{\Gamma
(2k)}{m!\Gamma (m+2k)},
\end{equation}
and to the Heisenberg-Weyl group
\begin{equation}\label{ort2}
\varphi_n=(n!)^{-\frac{1}{2}}(a^+)^{n}\varphi_0;~~<\varphi_{n'},\varphi_n>=\delta_{nn'}.
\end{equation}

We write down  the vector $e_0$ in (\ref{cond}) as 
\begin{equation}\label{zeroM}
 e_0=e_0^{H}\otimes e_0^{K'},~ \text{where}~\varphi_0\equiv
  e^H_0;~ e_0^{K'}\equiv e_{k',k'} .
\end{equation}
\begin{Proposition}\label{cevan}

The kernel
$K(z,w;\bar{z}',\bar{w}'):=(e_{\bar{z},\bar{w}},e_{\bar{z}',\bar{w}'}):
\mc{D}^J_1\times \bar{\mc{D}}^J_1\rightarrow \C$ is: 
\begin{equation}\label{KHK}
K(z,w;\bar{z}',\bar{w}')
=(1-{w}\bar{w}')^{-2k} \exp (F (z,w;\bar{z}',\bar{w}'));~ F= {\frac{2\bar{z}'{z}+z^2\bar{w}'
+\bar{z}'^2w}{2(1-{w}\bar{w}')}} ,
\end{equation}
\begin{equation}\label{hot}
K = K(\bar{z},\bar{w},z,w) =
(1-w\bar{w})^{-2k}\exp{\frac{2z\bar{z}+z^2\bar{w}+\bar{z}^2w}{2(1-w\bar{w})}},
~ z ,w\in\C, ~|w|<1 . 
\end{equation}

The normalized squeezed state vector   and the
un-normalized 
 Perelomov's coherent state vector
are related by the relation
\begin{equation}\label{csv}
\Psi_{\alpha, w} = (1-w\bar{w})^k
\exp (-\frac{\bar{\alpha}}{2}z)e_{z,w},~ z=\alpha-w\bar{\alpha}.
\end{equation}

The composition law in the Jacobi group 
$ G^J_1:=HW\rtimes SU(1,1 )$ is 
\begin{equation}\label{compositie}
(g_1,\alpha_1,t_1)\circ (g_2,\alpha_2, t_2)= (g_1\circ g_2,
g_2^{-1}\cdot\alpha_1+\alpha_2, t_1+ t_2 +\Im
(g^{-1}_2\cdot\alpha_1\bar{\alpha}_2)),
\end{equation}
where $g_{i}, i=1,2$ are of the form \begin{equation}\label{ggg}
g=\left(\begin{array}{cc}a & b \\\bar{b} & \bar{a}\end{array}\right),
|a|^2-|b|^2=1, 
\end{equation}
 $g\cdot\alpha :=\alpha_g$ is given by 
$\alpha_g = a\,\alpha + b\,\bar{\alpha}$, 
$g^{-1}\cdot\alpha =\bar{a}\alpha -b\bar{\alpha}$. 

Let $ (g,\alpha )\in G^J_1$ and let $(z,w)\in
{\mc{D}^J_1}:=\C\times\mc{D}_1$. The action of the group $G^J_1$ on the
manifold $\mc{D}^J_1$ is given by
\begin{equation}\label{xxx}
z_1=\frac{\alpha-\bar{\alpha}w+z}{\bar{b}w+\bar{a}}; ~ w_1=g\cdot w
=\frac{a w+ b}{\bar{b}w+\bar{a}}, ~g=\left(\begin{array}{cc} a& b\\c &d\end{array}\right)\in{\mbox{\emph{SU}}}(1,1).
\end{equation}
The scalar product of
functions from the space $\mc{F}_K$ corresponding to the kernel
defined by \emph{(\ref{KHK})}  on the manifold $\mathcal{D}^J_1$, with
$G^J_1$-invariant measure $d\nu$, is
\emph{(}$f_{\psi}(z,w)=(e_{\bar{z},\bar{w}},\psi)_{\Hi}$\emph{):}
\begin{equation}\label{ofi}
(\phi ,\psi )= \Lambda_1\! \int_{z\in\C;
|w|<1}\!\bar{f}_{\phi}(z,w)f_{\psi}(z,w)
(1\!-\!w\bar{w})^{2k}\!\exp{\!-\frac{|z|^2}{1\!-\!w\bar{w}}}
\exp{\!-\frac{z^2\bar{w}\!+\!\bar{z}^2w}{2(1\!-\!w\bar{w})}} d \nu ,
\end{equation}
\begin{equation}\label{ofi3}
d \nu =\frac{d \Re w d\Im w}{(1-w\bar{w})^3}d \Re z d \Im z, ~\Lambda_1 = \frac{4k-3}{2\pi^2} .
\end{equation}

The K\"ahler potential $f:=\log K$  is 
\begin{equation}\label{keler}
f =\frac{2z\bar{z}+z^2\bar{w}+\bar{z}^2w}{2(1-w\bar{w})}
-2k\log (1-w\bar{w}),
\end{equation}
and the  the K\"ahler two-form $\omega_1$ on $\mc{D}^J_1$, $G^J_1$-invariant to the
action \emph{(\ref{xxx})},  is 
\begin{equation}\label{aab}
-\i\: \omega_1 =\frac{2k}{(1-w\bar{w})^2}dw \wedge d\bar{w} +
\frac{A\wedge \bar{A}}{1-w\bar{w}},
~A=dz+\bar{\eta} dw, ~\eta=\frac{z+\bar{z}w}{1-w\bar{w}}.
\end{equation}
\end{Proposition}

\subsection{The real Jacobi group}\label{REALGR}
Let us now  recall that $C^{-1}\text{SL}_2(\R )C=\text{SU}(1,1)$,
where 
\begin{equation}\label{xcs}
C=\left(\begin{array}{cc} \i& \i \\
-1& 1\end{array}\right); C^{-1}= \frac{1}{2\i}\left(\begin{array}{cc}1
& -\i \\ 1 & \i \end{array}\right) . 
\end{equation}
If $M\in \SL$ is the matrix 
\begin{equation}\label{mst}
M=\left(\begin{array}{cc}a &\ b \\ c  &
d  \end{array}\right),~ad-bc = 1, 
\end{equation}
 then 
\begin{equation}\label{mstar}
M_*=C^{-1}MC=\left(\begin{array}{cc}\alpha&\beta\\ \bar{\beta} &
\bar{\alpha}\end{array}\right),~
\alpha,\beta\in\C, |\alpha|^2-|\beta|^2=1.
\end{equation}

The map (\ref{xcs}) induces a transformation of the bounded
domain $\mc{D}_1$ into the upper half plane $\mc{X}_1$ and
\begin{equation}\label{taustar}
w=C^{-1}(v )=\frac{v-\i}{v+\i}\in\mc{D}_1;~ v= Cw= \i\frac{1+w}{1-w}.
 \end{equation}

  K\"ahler and Berndt have investigated the
Jacobi group  $G^J_1(\R ):= \text{SL}_2(\R )\ltimes \R^2$ acting on the
Siegel-Jacobi upper half plane  $\mc{X}^J_1:=\mc{X}_1\times\C$
\cite{bern,bb,bs,cal1,cal2,cal3,cal},  where
$\mc{X}_1$ is the Siegel upper half plane
$\mc{X}_1:=\{v\in\C|\Im (v)>0\}$.  It is easy to proof (see also
(\cite{jac1}) the following


 \begin{Remark}\label{actnea1}
The action $C^{-1}G^J_1(\R )C$ 
 descends on the basis to the biholomorphic map: 
$\check{C}^{-1}: \mc{X}^J_1:= \mc{X}_1\times \C \rightarrow
\mc{D}^J_1:=\mc{D}_1\times \C$:
\begin{equation}\label{noiec}
w=\frac{v-\i}{v+\i};~ z=\frac{2\i u}{v+\i},~
w\in\mc{D}_1,~v\in\mc{X}_1,z\in\C .
\end{equation}
Under the partial Cayley transform \emph{(\ref{noiec})}, the K\"ahler  two-form
$\omega_1$  \emph{(\ref{aab})}
 becomes 
\begin{equation}\label{ura3}
-\i\; \omega_1' = -\frac{2k}{(\bar{v}-v)^2}dv\wedge
d\bar{v}+\frac{2}{\i(\bar{v}-v)}B\wedge \bar{B},  ~B=
du-\frac{u-\bar{u}}{v-\bar{v}}dv.
\end{equation}
$\omega_1'$ is K\"ahler homogeneous under the action of $G^J_1(\R)$ on
$\mc{X}^J_1$, $((h,\alpha),(v,u))\rightarrow (v_1,u_1)$:  
\begin{equation}\label{avc}
v_1=h\cdot v=\frac{av+b}{cv+d}, ~u_1=\frac{u+nv+m}{cv+d},
~h=\left(\begin{array}{cc}a & b\\ c&
    d\end{array}\right)\in\emph{\text{SL}}(2,\R),~\alpha=m+\i n,
\end{equation}
where the matrices $g$ in \emph{(\ref{xxx})} and $h$ in \emph{(\ref{avc})} are
related by \emph{(\ref{mstar})}.




\end{Remark}
\section{Fundamental
  Conjecture for the Siegel-Jacobi  disk and domain} \label{ok}
Firstly, we fix the terminology \cite{bo}.  A complex analytic manifold is
{\it K\"ahlerian}  if it is endowed with a Hermitian metric whose imaginary
part $\omega$  has $d\omega = 0$. A coset space is {\it homogenous
K\"ahlerian} if it caries a K\"ahlerian structure invariant under the
group. We call {\it a homogeneous K\"ahler diffeomorphism}  a
diffeomorphism $\phi:M\rightarrow N$  of homogeneous
K\"ahler manifolds  such that $\phi^*\omega_N=\omega_M$.

Let us remind the {\it fundamental conjecture for homogeneous K\"ahler
manifolds}
 (Gin\-dikin -Vinberg): {\it  every homogenous K\"ahler manifold is a holomorphic fiber bundle over a homogenous bounded domain in which the fiber is the product of a flat homogenous K\"ahler manifold and a compact simply connected homogenous K\"ahler manifold}. The compact case was considered by Wang
\cite{wa}, Borel \cite{bo} and   Matsushima \cite{ma} have  considered the
case of a transitive reductive group  of automorphisms,  while Gindikin and Vinberg  \cite{GV}
considered the transitive  automorphism group. We mention also the
essential contribution of Piatetski-Shapiro in this field
\cite{pia}. The complex version, in the formulation of Dorfmeister and
Nakajima \cite{DN}, essentially asserts that: {\it every homogenous
  K\"ahler manifold, as a complex manifold, is the product of a
  compact simply connected homogenous manifold (generalized flag
  manifold), a homogenous bounded domain,  and $\C^n/\Gamma$, where $\Gamma$ denotes a discrete subgroup of translations of $\C^n$}.

\begin{Proposition}\label{MBNM}
Let us  consider 
the  K\"ahler
two-form $\omega_1$  \emph{(\ref{aab})},  $G^J_1$-invariant invariant under the action
\emph{(\ref{xxx})} of $G^J_1$ on the homogenous K\"ahler Siegel-Jacobi disk  $\mc{D}^J_1$. We have the
homogenous K\"ahler diffeomorphism  $$FC:~ (\mc{D}^J_1,\omega_1)\rightarrow 
(\mc{D}_1\times\C, \omega_0)=
(\mc{D}_1,\omega_{\mc{D}_1} )\otimes
(\C,\omega_{\C}), ~\omega_0=FC^*(\omega_1),$$
\begin{equation}\label{zZx} FC: ~z=\eta-w \bar{\eta},~
 FC^{-1}: ~ \eta=\frac{z+w\bar{z}}{1-|w|^2},
\end{equation}
\begin{equation}\label{mtr0}
\omega_0= \omega_{\mc{D}_1} + \omega_{\C};~ -\i \omega_{\mc{D}_1}= 
\frac{2k}{(1-w\bar{w})^2}dw\wedge d\bar{w}, ~ -\i
\omega_{\C}=d\eta\wedge d\bar{\eta}.
\end{equation}  

The K\"ahler two-form \emph{(\ref{mtr0})} is invariant  at the action of
$G^J_1$ on $\C\times\mc{D}_1$, $((g,\alpha), (\eta,w))\rightarrow (\eta_1,w_1)$,
\begin{equation}\label{grozav}
\eta_1=a(\eta+\alpha)+b(\bar{\eta}+\bar{\alpha}), ~
w_1=\frac{aw+b}{\bar{b}w+\bar{a}}, ~g=\left(\begin{array}{cc} a & b\\
      c& d\end{array}\right)\in {\emph{\text{SU}}}(1,1). 
\end{equation}

We have also the homogenous K\"ahler diffeomorphism 
$$FC_1:(\mc{X}^J_1,\omega'_1)\rightarrow (\mc{X}\times\C, \omega'_0) =
(\mc{X}_1, \omega_{\mc{X}_1})\times(\C,\omega_{\C}),
~\omega'_0=FC_1^*(\omega'_1) , $$  
\begin{equation}\label{XCX}
FC_1: 2\i u =(v+\i)\eta-(\bar{v}-\i)\bar{\eta}; ~~FC^{-1}_1:
\eta=\frac{u\bar{v}-\bar{u}v+\i (\bar{u}-u)}{\bar{v}-v}, 
\end{equation}
where $\omega'_1$  is the K\"ahler two-form \emph{\text{(\ref{ura3})}},
$G^J_1(\R)$-invariant to the action \emph{(\ref{avc})},    and
\begin{equation}\label{mtr0p}
  \omega'_0= \omega_{\mc{X}_1}+\omega_{\C}, ~ \i  d\omega_{\mc{X}_1}=
  \frac{2k}{(v-\bar{v})^2}dv\wedge d\bar{v}. 
\end{equation} 
\end{Proposition}

{\it Proof}. 
The idea is to use the transformation (\ref{noiec}) and the EZ
(Eichler-Zagier) coordinates (\ref{eqez}), 
(cf. the definition at p. 12 and p. 51
in \cite{bs}
adapted to our notation)
\begin{equation}\label{eqez}
v=x+iy; ~ u=pv+q,~ x,p,q, y\in\R ,y>0,
\end{equation}
and come back from $v$ to $w$. So,
let $$z=2\i\frac{u}{v+\i}=2\i\frac{pv+q}{v+\i},$$  where, by the
(inverse Cayley) transform (\ref{taustar}),  $v=-\i\frac{w+1}{w-1}$. We have $z=q+\i p +w(-q+\i p)$, and if denote $\eta = q+\i p$, where $q,p\in\R$, then 
$z=\eta-w\bar{\eta}$, with $\eta$  appearing already in  (\ref{aab}),
and $A=d\eta-w d\bar{\eta}$. The last term in (\ref{aab}) becomes
\begin{equation}\frac{A\wedge\bar{A}}{1-|w|^2}=d\eta\wedge d\bar{\eta}=2\i dp\wedge dq.\end{equation}
Vice-versa, we have $
 d\eta =\frac{A+w\bar{A}}{1-|w|^2}$,
with $A$ given in (\ref{aab}). 

For the second assertion, we introduce 
 the transformation (\ref{noiec})     $ z=2\i
u(v+\i)^{-1}$  in (\ref{zZx})  and we get:
$2\i (u-\bar{u}) =(\eta-\bar{\eta})(v-\bar{v})$. 
Than  $B$ in (\ref{ura3})  becomes
$$B=\frac{1}{2\i}[(v+\i)d\eta-(v-\i)d\bar{\eta}]$$
and we get (\ref{mtr0p}). \hfill $\gata$

\begin{corollary}\label{siasta}
Let us denote $\mc{F}:=F\circ FC$, $\mc{K}=K\circ FC$. In the  variables $(\eta, w)$, the scalar product
 \emph{\text{ (\ref{KHK})}}
becomes
\begin{equation}\label{HOTt}
K(w,\eta;\bar{w}',\bar{\eta}')=(1-w\bar{w}')^{-2k}\exp \mathcal{F},\quad\emph{where}
\end{equation}
$$
2\mathcal{F}=2(\bar{\eta}\zeta+|\eta'|^2)-w\bar{\eta}^2-\bar{w}'\eta'^2+
(1-w\bar{w}')^{-1}(-2|\zeta|^2+w{\bar{\zeta}}^2+\bar{w}'\zeta^2), ~\zeta=\eta-\eta'$$
and, for $w=w',\zeta=0$, \emph{\text{ (\ref{hot})}} becomes
\begin{equation}\label{hott}
\mc{K}=(1-w\bar{w})^{-2k}\exp \mathcal{F}, ~ 2 \mathcal{F}=  2 \eta\bar{\eta}-\bar{w}\eta^2-w\bar{\eta}^2.
\end{equation}
Under the $FC$ transform, the scalar product  \emph{\text{(\ref{ofi})}} on $\mc{D}^J_1$, with $G^J_1$-invariant
measure $d\nu$ at the action \emph{(\ref{xxx})}, becomes the scalar product  \emph{(\ref{ofiXXXX})} on
$\C\times \mc{D}_1$, with  the measure $d\nu'$
\emph{(\ref{ofi333})},  $G^J_1$-invariant
at the action \emph{(\ref{grozav})}
 (also $d\omega_0\wedge d\omega_0= -8kd\nu'$): 
\begin{equation}\label{ofiXXXX}
(\phi ,\psi )= \Lambda_1\! \int_{\eta\in\C;
|w|<1}\!\bar{f}_{\phi}(\eta,w)f_{\psi}(\eta,w)
(1\!-\!w\bar{w})^{2k}\!\exp(-\mathcal{F}) d \nu' , \mbox{\quad\emph{where}}
\end{equation}
\begin{equation}\label{ofi333}
d \nu' =\frac{d \Re w d\Im w}{(1-w\bar{w})^2}d \Re \eta d \Im \eta .
\end{equation}
 \end{corollary}                                                                                                                                     

\subsection{Geodesics on $\mathcal{D}^J_1$}\label{GODZ}
Now we look at the effect of the $FC$ transform on the equations of
geodesics on  $\mathcal{D}^J_1$. We recall (cf. \cite{jac1}) 
\begin{Remark}
The equations of the geodesics on the manifold $\mc{D}^J_1$, 
endowed with the two-form {\em({\ref{aab}})} in the variables $(w,z)\in\mc{D}_1\times\C$,
are
\begin{subequations}
\begin{eqnarray}\label{geod}
2k\frac{d^2z}{dt^2}-\bar{\eta}(\frac{dz}{dt})^2+2(2k\frac{\bar{w}}{P}-
\bar{\eta}^2)\frac{dz}{dt}\frac{dw}{dt}-\bar{\eta}^3(\frac{dw}{dt})^2
&= & 0;\label{eq1g}\\
2k\frac{d^2w}{dt^2}+(\frac{dz}{dt})^2+2\bar{\eta}\frac{dz}{dt}\frac{dw}{dt}
+ (4k\frac{\bar{w}}{P}+\bar{\eta}^2)(\frac{dw}{dt})^2 & = & 0, \label{eq2g}
\end{eqnarray}
\end{subequations}
where $\eta$ is given by  {\em{(\ref{zZx})}} and
$P=1-w\bar{w}$.\end{Remark} 
If we introduce  the
solution  \begin{equation}\label{geow}
w=w(t)=B\frac{\tanh(t\sqrt{\bar{B}B})}{\sqrt{\bar{B}B}}\end{equation}
of the equations  of geodesics on $\mathcal{D}_1$,  
$$\frac{d^2w}{dt^2}+2\frac{\bar{w}}{1-w\bar{w}}(\frac{dw}{dt})^2 =0,  $$
into  (\ref{eq2g}), then (\ref{eq2g}) becomes 
\begin{equation}\label{1primg}
(\frac{dz}{dt}+\bar{\eta}\frac{dw}{dt})^2=0.
\end{equation}
Now we introduce the solution $\frac{dz}{dt}=-\bar{\eta}\frac{dw}{dt}$ of  (\ref{1primg}) into (\ref{eq1g}) and we obtain:
\begin{equation}\label{2primg}
\frac{d^2z}{dt^2}-\frac{2\bar{w}}{P}\bar{\eta}(\frac{dw}{dt})^2=0. 
\end{equation}
If in (\ref{2primg}) we take into account (\ref{1primg}), we get
$\frac{d\bar{\eta}}{dt}=0$, and  {\it a particular solution of} (\ref{geod})
{\it consists of} $(\eta = ct,  w)$  {\it with} $w$ {\it given by}
(\ref{geow}).
\hfill $\gata$

\section{Classical motion and quantum evolution}\label{CLSQ}

Let $M=G/H$ be a homogeneous  manifold with a $G$-invariant K\"ahler two-form
$\omega$
\begin{equation}\label{kall}
\omega(z)=\i\sum_{\alpha\in\Delta_+} g_{\alpha,\beta}  d z_{\alpha}\wedge
d\bar{z}_{\beta}, ~g_{\alpha,\beta}=\frac{\pa^2}{\pa
  z_{\alpha} \pa\bar{z}_{\beta}} \log <e_z,e_z>. 
\end{equation}
Above $e_z\in\Hi$ are Perelomov's coherent state vectors,  indexed by the
points $z\in M$,  obtained by  the unitary irreducible representation
$\pi$  on $\Hi$ of
$G$,  $e_z=\exp(\sum_{\alpha\in\Delta_+}z_{\alpha}X_{\alpha})$, and 
$\Delta_+$ are the positive roots of the Lie algebra $\got{g}$ of $G$,
with generators $X_\alpha,\alpha\in\Delta$ \cite{perG}. 

Passing on from the dynamical system problem
 in the Hilbert space $\Hi$ to the corresponding one on $M$ is called
sometimes {\it dequantization}, and the dynamical system on $M$ is a classical
one \cite{sbcag,sbl}. Following Berezin \cite{berezin2},\cite{berezin1}, the
motion on the classical phase space can be described by the local
equations of motion
$\dot{z}_{\alpha}=\i \left\{\mc{H},z_{\alpha}\right\},
  ~\alpha \in \Delta_+ $,
where $\mc{H}$ is
  the classical Hamiltonian $\mc{H}=<e_z,e_z>^{-1}<e_z|\mb{H}|\e_z>$
  (the covariant 
  symbol) attached to
  the quantum Hamiltonian $\mb{H}$, and the Poisson bracket is
  introduced using the matrix $g^{-1}$. 

We consider an algebraic Hamiltonian linear in the generators
${\mb{X}}_{\lambda} $ of the
group of symmetry $G$
\begin{equation}\label{lllu}
\mb{H}=\sum_{\lambda\in\Delta}\epsilon_{\lambda}{\mb{X}}_{\lambda} .
\end{equation}
The  classical motion generated
by the Hamiltonian  (\ref{lllu}) is  given by the  equations
of motion  on $M=G/H$ \cite{sbcag,sbl}: 
\begin{equation}\label{moveM}
{\i\dot{z}_{\alpha}=\sum_{\lambda\in\Delta}\epsilon_{\lambda}Q_{\lambda
,\alpha}},~\alpha\in\Delta_+ , 
\end{equation}
where the differential action corresponding to the operator
$\mb{X}_{\lambda}$ in (\ref{lllu}) can be expressed in a local
system of coordinates as a holomorphic  first order differential
operator with polynomial coefficients  ($\pa_{\beta}=\frac{\pa}{\pa
  z_{\beta}}$),
\begin{equation}\label{VBC}\db{X}_{\lambda}=P_{\lambda}+\sum_{\beta\in\Delta_+}Q_{\lambda,
    \beta}\partial_{\beta}, \lambda\in\Delta.
\end{equation}

We look also for the solutions of the Schr\"odinger equations attached
to the Hamiltonian  $\mb{H}$  (\ref{lllu}) 
\begin{equation}\label{SCH}
\mb{H}\psi = \i \dot{\psi},~~
\quad{\mbox{where~~~}}\psi= \e ^{\i \varphi}<e_z,e_z>^{-1/2}e_z.
\end{equation}

We remember that  \cite{sbl}
\begin{Proposition}\label{MERETH}
On the homogenous  manifold $M = G/H $ on which the holomorphic representation \emph{\text{(\ref{VBC})}} is true, the classical
motion and the quantum evolution generated by the linear Hamiltonian
\emph{\text{(\ref{lllu})}}  are given by the same equation of motion
\emph{\text{(\ref{moveM})}}. The phase $\varphi$ in \emph{\text{(\ref{SCH})}} is given by the
sum $\varphi=\varphi_D +\varphi_ B$ of the
dynamical and Berry phase, 
\begin{subequations}
\begin{eqnarray}
\varphi_D & =& -\int_0^t\mc{H}(t)dt, \quad{\mbox{\emph{where}~~~}}\label{phiD}\\
\mc{H}(t) & = &\sum_{\lambda\in\Delta}\epsilon_{\lambda}(P_{\lambda}+
 \sum_{\beta\in\Delta_+}Q_{\lambda,\beta}\pa_{\beta}\ln <e_z,e_z>) \nonumber\\
 & = &
\sum_{\lambda\in\Delta}\epsilon_{\lambda}P_{\lambda}+
\i\sum_{\beta\in\Delta_+}\dot{z}_{\beta}\pa_{\beta}\ln <e_z,e_z>,\quad{\mbox{~~~}}\nonumber\\
\varphi_B & = &   -\Im\int_0^t<e_z,e_z>^{-1}<e_z|d|e_z>\label{phiB}\\
 & = & \nonumber
 \frac{\i}{2}\int_0^t\sum_{\alpha\in\Delta_+}(\dot{z}_{\alpha}\pa_{\alpha}-\dot{\bar{z}}_{\alpha}\bar{\pa}_{\alpha} ) \ln <e_z,e_z>.
\end{eqnarray}
\end{subequations}
\end{Proposition}

\subsection{Equations of motion on Siegel-Jacobi disk and domain}
Let us consider a linear hermitian Hamiltonian in the generators of the Jacobi
group $G^J_1$:
\begin{equation}\label{guru}
\mb{H} = \epsilon_a\mb{a} +\bar{\epsilon}_a\mb{a}^{\dagger}
 +\epsilon_0 {\mb{K}}_0 +\epsilon_+{\mb{K}}_++\epsilon_-{\mb{K}}_-  ,~~
\bar{\epsilon}_+=\epsilon_-, ~ \epsilon_0=\bar{\epsilon_0}.
\end{equation}

With Lemma \ref{mixt},  Proposition \ref{MERETH} and (\ref{hot}), we
get 
\begin{Proposition}\label{EQLIN1}
The equations of motion on the Siegel-Jacobi disk $\mathcal{D}^J_1$
 generated by
the linear Hamiltonian {\em(\ref{guru})}
 are:
\begin{subequations}\label{qqqN}
\begin{eqnarray}
i\dot{z} & = & \epsilon_a+{\bar{\epsilon}}_a w+(\frac{\epsilon_0}{2}
+\epsilon_+  w )z,~z,w\in\C, ~|w|<1, \label{guru2}\\
i\dot{w} & = & \epsilon_- + \epsilon_0w+
\epsilon_+w^2 .\label{guru1}
\end{eqnarray}
\end{subequations}
The equations of motion generated by the linear Hamiltonian
{\emph{(\ref{guru})}} on the manifold $\mc{X}^J_1$, obtained from the
equations  \emph{(\ref{qqqN})} by partial Cayley transform  \emph{(\ref{noiec})} are
\begin{subequations}\label{qqq1}
\begin{eqnarray}
~~~- 2\dot{v} & = &
(\epsilon_0+\epsilon_++\epsilon_-) v^2+2\i
(\epsilon_--\epsilon_+)v +\epsilon_0-\epsilon_--\epsilon_+,
~v\in\C,\Im v>0,
\label{drac}\\
~~~-2\dot{u} & = & (\epsilon_a+\bar{\epsilon}_a)v +\i(\epsilon_a-\bar{\epsilon}_a)
+[(\epsilon_0+\epsilon_++\epsilon_-)v +\i (\epsilon_--\epsilon_+)]u,~ u\in\C .\label{dracu1}
\end{eqnarray}
\end{subequations}

 For the $\eta$ defined in the $FC^{-1}$ transform \emph{(\ref{zZx})}, the
 system of first order differential  equations \emph{(\ref{qqqN})} becomes  the
system of separate equations 
\begin{subequations}\label{qqqN1}
\begin{eqnarray}
i\dot{\eta} & = & \epsilon_a
+\epsilon_-\bar{\eta}+\frac{\epsilon_0}{2}\eta,~\eta \in\C,\label{guru22}\\
i\dot{w} & = & \epsilon_- + \epsilon_0w+
\epsilon_+w^2,~ w\in\C,|w|<1, .\label{guru12}
\end{eqnarray}
\end{subequations}

 If in \emph{(\ref{dracu1})}  we make the
change of variables \emph{(\ref{XCX})},  we get the system of decoupled equations of motion
on $\mc{X}^J_1$
\begin{subequations}\label{qqqN11}
\begin{eqnarray}
i\dot{\eta} & = & \epsilon_a
+\epsilon_-\bar{\eta}+\frac{\epsilon_0}{2}\eta,~\eta\in\C,  \label{guru2222}\\
- 2\dot{v} & = &
(\epsilon_0+\epsilon_++\epsilon_-) v^2+2\i
(\epsilon_--\epsilon_+)v +\epsilon_0-\epsilon_--\epsilon_+, ~v\in\C, \Im v>0,
\end{eqnarray}
\end{subequations}
\end{Proposition}
The  equation (\ref{guru1})  (the  equation
 (\ref{drac})) 
 is a {\it Riccati equation} on $\mc{D}_1$ (respectively,  on $\mc{X}_1$).  Remark that the dynamics
on the Siegel ball $\mc{D}_1$,  determined by the Hamiltonian
(\ref{guru}),  {\it linear in the generators of the Jacobi group
  $G^J_1$, 
depends {\bf only} on the generators of the group}  $\text{SU}(1,1)$. The
Riccati equation  on the $\mc{D}_1$
(\ref{guru1}) appears in literature, see e.g. equation (18.2.8) in
\cite{perG} in the context of quantum oscillator with variable frequency.


{\bf 1.a} 
We consider the case of {\it constant coefficients} of the Hamiltonian  (\ref{guru}). There are two equivalent methods to integrate the Riccati equation
(\ref{guru1}).

In the equation (\ref{guru1}) we put 
$w=-\frac{\i}{\epsilon_+}\frac{\dot{a}}{a}$, 
and  we get for $a$ the equation 
$\ddot{a}+\i\epsilon_0\dot{a}-\epsilon_-\epsilon_+a=0$,
which has the characteristic equation  ($a(t)=C_{1,2}\e^{\i w_{1,2}t}$)
\begin{equation}\label{CHAR}
w^2_{1,2}+\epsilon_0w_{1,2}+\epsilon_+\epsilon_-=0;
w_{1,2}=\frac{-\epsilon_0\pm\sqrt{\Delta}}{2},
\Delta=\epsilon_0^2-4\epsilon_+\epsilon_- .
\end{equation}

For $\epsilon$-s  constant in (\ref{guru1}), the solution of the
Riccati equation is
\begin{equation}\label{WcomP}
w(t) = \frac{1}{\epsilon_+}\cdot\frac{w_1C_1\e^{\i w_1t}+w_2C_2\e^{\i
    w_2t}}{C_1\e^{\i w_1t}+C_2\e^{\i w_2t}}= 
 \frac{1}{\epsilon_+}\cdot\frac{w_1C_1\e^{\frac{\i\sqrt{\Delta}}{2}t}+w_2C_2\e^{{-\frac{\i\sqrt{\Delta}}{2}t}}}{C_1\e^{\frac{\i\sqrt{\Delta}}{2}t}+C_2\e^{{-\frac{\i\sqrt{\Delta}}{2}t}}},
\end{equation}
where $w_{1,2}$ are given in  (\ref{CHAR}), and
in order to have $w\in\C$, we must have 
$\Delta > 0 $  so, if also $\epsilon_0>0$, then  $w_{1,2}<0$. 
Imposing to the solution  of the Riccati equation
(\ref{guru1}) the initial condition $w(0)=w_0$, it results for
$f=\frac{C_1}{C_2}$ the value
$f =\frac{\epsilon_+w_0-w_2}{w_1-\epsilon_+w_0}$, 
and  we rewrite down  the solution (\ref{WcomP}) of (\ref{guru1}) as 
\begin{equation}\label{sfinal1}
w(t,w_0)=\frac{1}{\epsilon_+} \cdot\frac{fw_1\e^{\i\sqrt{\Delta}t}+w_2}{1+f\e^{\i\sqrt{\Delta}t}}. 
\end{equation}
For the solution (\ref{WcomP}) of the differential equation with
constant coefficients we find out   
$$1-w\bar{w}=\frac{\sqrt{\Delta}}{\epsilon_+\epsilon_-}\frac{-w_1|C_1|^2+
w_2|C_2|^2}
{|C_1\e^{\i\frac{\sqrt{\Delta}}{2}t}+C_2\e^{-\i\frac{\sqrt{\Delta}}{2}t}|^2},$$
and {\it the condition} $w(t)\in\mc{D}_1$ 
 {\it imposes the restrictions}:
\begin{equation}\label{restict}
\vert\frac{C_1}{C_2}\vert>\sqrt{\frac{w_2}{w_1}}= 
\frac{1+\sqrt{1-\delta}}{\sqrt{\delta}},~\epsilon_0>0, ~\Delta >0, 
~\delta =4\frac{\epsilon_+\epsilon_-}{\epsilon_0^2} <1. 
\end{equation}
The second method to integrate the Riccati equation is to
make the substitution $w=X/Y$, and we associate to (\ref{guru1}) the linear system  of first
order differential equations
\begin{equation}\label{EQUIV}
\i \dot{X}  ~ =  ~\epsilon_-Y+\epsilon_0X, ~~~
\i \dot{Y}  ~  =  ~  -\epsilon_+X . 
\end{equation}
We eliminate $X$,  integrate the equation in $Y$,  and
finally, we get for $w$ the same solution (\ref{WcomP}).

{\bf 1.b}  Integration of equation (\ref{drac}). 

We write (\ref{drac}) as 
\begin{equation}\label{dracmic}
-\dot{v}=Av^2+Bv+C,  
~ A=\frac{1}{2}(\epsilon_-+\epsilon_++\epsilon_0);~ 
B=\i(\epsilon_--\epsilon_+); ~
C=\frac{1}{2}(\epsilon_0-\epsilon_--\epsilon_+),
\end{equation}
 where $A,B,C\in\R$. 

We look for a
solution of (\ref{dracmic})  as $v=\frac{1}{A}\frac{\dot{b}}{b}$  and we
have 
$\ddot{b}+B\dot{b}+ACb=0$,
which has the solution $b(t)=C'_{1,2}\e^{\i v_{1,2}t}$, with $v_{1,2}=\frac{\epsilon_+-\epsilon_-\pm\sqrt{\Delta}}{2},
  ~\Delta=\epsilon_0^2-4\epsilon_+\epsilon_- >0$,
The solution of (\ref{dracmic}) is 
\begin{equation}\label{VcomP}
v(t)=\frac{\i}{A}\cdot\frac{v_1C_1'\e^{\i v_1t}+v_2C_2'\e^{\i v_2t}}{C_1'\e^{\i
    v_1t}+C_2'\e^{\i
    v_2t}}=\frac{\i}{A}\cdot\frac{v_1C_1'\e^{\frac{\i\sqrt{\Delta}}{2}t}+
v_2C_2'\e^{-\frac{\i\sqrt{\Delta}}{2}t}}
{C_1'\e^{\frac{\i\sqrt{\Delta}}{2}t}+C_2'\e^{-\frac{\i\sqrt{\Delta}}{2}t}}
\end{equation}
which is complex in the same case as the solution (\ref{WcomP}),
$\Delta>0$.

Note that {\it the solution \emph{(\ref{VcomP})}  of the
Riccati equation
\emph{(\ref{drac})}  on $\mc{X}_1$ is  related to the solution
\emph{(\ref{WcomP})} of the Riccati equation \emph{(\ref{guru1})}  on $\mc{D}_1$  by
the Cayley transform \emph{(\ref{taustar})}  if we chose the constants such
that}
\begin{equation}\label{CMNV}
\frac{C_1}{C_2}\cdot\frac{C'_2}{C'_1}=\frac{\epsilon_+-w_2}{\epsilon_+-w_1}.
\end{equation}

{\bf 2.} Suppose that we know the general solution of the Riccati
equation (\ref{guru1}) obtained integrating the {\it time-dependent} linear
system (\ref{EQUIV}).  Then (\ref{guru2}) becomes
\begin{equation}\label{nonH}
\i \dot{z}(t)= A(t) + B(t)z,  ~\text{where}
                               ,~ A(t)=\epsilon_a+\bar{\epsilon}_aw, ~
                                B(t)=
                               ~\frac{\epsilon_0}{2}+\epsilon_+w.
\end{equation}
Let $z(t)=CF(t)$ be the solution of the homogeneous equation $\i
\dot{z}(t)= B(t)z$, where $F(t)=\exp (-\i\int_{t_0}^tB(t)d
t)$.  Then the solution of the  differential  equation (\ref{dracmic})  is
given by 
$C(t)=C_0-i\int_{t_0}^t\frac{A}{F}dt.$

{\bf 3.}  Now we look at the {\it  decoupled} system of differential
equations (\ref{qqqN1}), also in the {\it autonomous case}.
We write down the complex numbers as 
$\eta=x+\i  y$, $ \epsilon_a=
a+\i b$, $ \epsilon_-=m+\i n$ ,  $\epsilon_0/2=p$, 
and
(\ref{guru22}) becomes  the linear system  of differential equations
\begin{equation}\label{qqXY}
\dot{x}  ~ = ~  nx +(p-m)y+b, ~~
\dot{y} ~  = ~  -(m+p)x-ny-a , 
\end{equation}
where the solution of the characteristic equation in
 $\e^{(\lambda t)}$: 
\begin{equation}\label{CHARRE}
\det\left(\begin{array}{cc} n-\lambda & p-m\\ -(p+m) &
   -( n+\lambda) \end{array}\right) =0,~ \lambda^2 = n^2+m^2-p^2, 
\end{equation}
is $\lambda
=\pm\frac{\i}{2}\sqrt{\Delta}$. 

The solution of (\ref{qqXY}) in the case of constant coefficients  is
\begin{subequations}\label{etat}
\begin{eqnarray}\label{ssisxy}
x(t) ~ & = & \frac{q}{2\lambda}(\alpha\e^{\lambda t}-\beta\e^{-\lambda t}) -\frac{q}{\lambda^2},~q=nb+a(m-p), ~\alpha,\beta\in\C,\\
y(t)~ & = & \frac{q}{2\lambda}\cdot\frac{\alpha(\lambda -n)\e^{\lambda
  t}+\beta(\lambda+n)\e^{-\lambda t}}{p-m}+\frac{\frac{nq}{\lambda^2}-b}{p-m}, ~\lambda
= \i \frac{\sqrt{\Delta}}{2},
\end{eqnarray}
\end{subequations}
and {\it we find the solution $\eta(t) = x(t) +\i y(t)$ of the differential equation} (\ref{guru22})
 \begin{subequations}\label{etass}
\begin{eqnarray}
\mbox{\qquad}
~ \eta(t)  ~ & = & M \e^{\i\frac{\sqrt{\Delta}}{2}t}+N\e^{-\i\frac{\sqrt{\Delta}}{2}t} +P, \quad{\text{where}}\\
\mbox{\qquad} 
 ~ M  ~ & = &-\i \frac{q \alpha}{r\sqrt{\Delta}}(\epsilon_-+w_1);
~ N ~  =  \i\frac{q \beta}{r \sqrt{\Delta}}(\epsilon_-+w_2), \\
~ \frac{\alpha}{\bar{\beta}} ~ & = &
\frac{\epsilon_-(\epsilon_++w_2)}{w_2(\epsilon_-+w_1)}
=\frac{w_1(\epsilon_++w_2)}{\epsilon_+(\epsilon_-+w_1)}, ~\alpha=\i\frac{r}{q}(\eta(t=0)-P),\\
 \mbox{\qquad} 
~ P ~  &= & \frac{4\epsilon_-\bar{\epsilon}_a-2\epsilon_0\epsilon_a}{\Delta}
, ~r  = \frac{1}{2}(\epsilon_-+\epsilon_+-\epsilon_o), 
\\
\mbox{\qquad} 
~ q ~  &= &- \frac{\epsilon_0}{4}(\epsilon_a+\bar{\epsilon}_a)+\frac{1}{2}(\epsilon_a\epsilon_++\bar{\epsilon}_a\epsilon_-).
\end{eqnarray}
\end{subequations}
The solution of the system of differential equations (\ref{qqqN}) is
given by $z=\eta -w\bar{\eta}$, where  $\eta(t)$  has the
expression 
given by (\ref{etass}), while the solution $w(t)$  of (\ref{guru1}) is given
by (\ref{WcomP}).

{\bf 4.} Instead of the linear hermitian Hamiltonian (\ref{guru}), we could consider
 the {\it non-hermitian Hamiltonian}
\begin{equation}\label{Guru}
\mb{H} = \epsilon_a\mb{a} +\epsilon_b\mb{a}^{\dagger}
 +\epsilon_0 {\mb{K}}_0 +\epsilon_+{\mb{K}}_++\epsilon_-{\mb{K}}_-  .
\end{equation}
which leads to the equations of motion (\ref{qqqN}), where in
(\ref{guru2}), the term linear in $w$ should have as coefficient
$\epsilon_b$ instead of $\bar{\epsilon}_a$ . Then, with the
change of variables $FC$ given by  (\ref{zZx}), we get instead of
(\ref{guru22} ) the equation 
\begin{subequations}
\begin{eqnarray}\label{GURU22}
~~~\i\dot{\eta} & = & (1-\bar{w}w)^{-1}(R + S \eta +T \bar{\eta})\label{GHK}, {\mbox{\quad}}\text{where}~
~ R =   \epsilon_a +(\epsilon_b-\bar{\epsilon}_a )w-\bar{\epsilon}_b\bar{w}w,\\
~~~S &  = & \frac{\epsilon_0}{2}+ (\epsilon_+-\bar{\epsilon}_-) w-\frac{\bar{\epsilon}_0}{2}\bar{w}w,
~~~T  =  \epsilon_- + (\frac{\epsilon_0-\bar{\epsilon}_0}{2})w-\bar{\epsilon}_+\bar{w}w .
\end{eqnarray}
\end{subequations}
 {\it The equation} (\ref{GHK}) {\it get the form} (\ref{guru22})
{\it only  if the
Hamiltonian} (\ref{Guru}) {\it becomes the hermitian Hamiltonian}
(\ref{guru}), i.e.  $\epsilon_b=\bar{\epsilon}_a$,
$\bar{\epsilon}_0=\epsilon_0$, and $\epsilon_-=\bar{\epsilon}_+$.

\subsection{Berry  phase for $\mathcal{D}^J_1$}   
     
We calculate Berry phase  with (\ref{phiB}), which on $\mc{D}^J_1$ reads
$$\frac{2}{\i}d\varphi_B= (dw\frac{\pa }{\pa w}- d\bar{w}\frac{\pa }{\pa
  \bar{w}}+dz\frac{\pa }{\pa z}- d\bar{z}\frac{\pa }{\pa
  \bar{z}}) f, $$          
where $f$ is  the K\"ahler potential (\ref{keler}). 

We have 
$$ f_{\bar{z}} = \eta, ~~ f_{\bar{w}} = \frac{\eta^2}{2}+\frac{2kw}{1-w\bar{w}}.$$
$$\frac{2}{\i}d\varphi_B=
(\frac{\bar{\eta}^2}{2}+\frac{2k\bar{w}}{1-w\bar{w}}) dw-
(\frac{\eta^2}{2}+\frac{2kw}{1-w\bar{w}}) d\bar{w}+\bar{\eta}dz -
\eta d\bar{z}. $$  
But $z=\eta- w\bar{\eta}$, and {\it the Berry phase on} $\mathcal{D}^J_1$
{\it in the variables} $(w,\eta)$ {\it is}
\begin{equation}\label{FFV}
\frac{2}{\i}d\varphi_B=
(\frac{2k\bar{w}}{1-w\bar{w}}-\frac{\bar{\eta}^2}{2}) dw
+(\bar{\eta} + \bar{w}\eta)d \eta - cc . 
\end{equation}
\subsection{Dynamical  phase}
The dynamical phase is calculated with (\ref{phiD}). Firstly we
calculate 
{\it the energy function} attached to the Hamiltonian (\ref{guru}): 
\begin{equation}\label{realH}
\begin{split}
\mc{H}= & ~ k\epsilon_0+ \bar{\epsilon}_az
+\epsilon_+(2kw+\frac{z^2}{2})+
  [\epsilon_a+\bar{\epsilon}_aw+(\frac{\epsilon_0}{2}+ \epsilon_+w)z]\bar{\eta}\\
&
+(\epsilon_-+\epsilon_0w+\epsilon_+w^2)(\frac{\bar{\eta}^2}{2}+\frac{2k\bar{w}}{1-w\bar{w}}).\\
\end{split}
\end{equation}
Now we put into  evidence that the energy function attached to the
hermitian Hamiltonian (\ref{guru})  is real and write it as  $\mc{H}=
\mc{H}_{\eta}+\mc{H}_{w}$, where 
\begin{subequations}\label{realHH}
\begin{eqnarray}
~\mc{H}_{\eta}  ~ & = & ~  \bar{\epsilon}_a\eta+\epsilon_a\bar{\eta} + 
\frac{1}{2}(\epsilon_+\eta^2+\epsilon_-\bar{\eta}^2+\epsilon_0\eta\bar{\eta}),\label{realHH1} \\
~ \mc{H}_{w} ~ & = & ~  k\epsilon_0+ \frac{2k}{1-w\bar{w}}(\epsilon_+w+\epsilon_-\bar{w}+\epsilon_0w\bar{w}). \label{realHH2}
\end{eqnarray}
\end{subequations}
Then the solutions of  (\ref{qqqN1}) are introduced into the expression of
  the energy function (\ref{realHH}). In the case of
  {\it constant coefficients} we use the solution (\ref{WcomP})  ((\ref{etass})) and we  find for the
   function (\ref{realHH2}) (respectively, (\ref{realHH1}))
the {\it  independent of time}  expression
\begin{subequations}\label{WRONG}
\begin{eqnarray}
\mc{H}_w(w(t)) & = & k(\epsilon_0+2\frac{-w_1^2|C_1|^2+w_2^2|C_2|^2}{-w_1|C_1|^2+w_2|C_2|^2}) ,\\ 
\mc{H}_{\eta}(\eta(t)) & = & \frac{2}{\Delta}\cdot
(\epsilon_+\epsilon^2_a+\epsilon_-\bar{\epsilon}_a^2-\epsilon_0|\epsilon_a|^2)-\frac{q^2}{r}|\alpha|^2 .
\end{eqnarray}
\end{subequations}

We look for the critical points of the energy function
(\ref{realHH}), i.e. the points $(w,\eta)\in\mc{D}_1\times\C$ for which
$\frac{\pa \mc{H}}{ \pa w} =0$, $\frac{\pa \mc{H}}{\pa \eta} = 0$: 
\begin{equation}\label{eqmis}
\frac{\pa \mc{H}}{\pa \eta}  = 
\bar{\epsilon_a}+\epsilon_+\eta+\frac{\epsilon_o}{2}\bar{\eta}= 0; ~~
\frac{\pa \mc{H}}{\pa  w} =  \frac{2k}{(1-w\bar{w})^2}(\epsilon_-\bar{w}^2+\epsilon_0\bar{w}+\epsilon_+)
  =  0. 
\end{equation}
The solution $w_c$ of (\ref{eqmis}),  also the solution of
$\dot{w}=0$, 
is 
 $ w_c=\frac{-\epsilon_0\pm\sqrt{\Delta}}{2\epsilon_+}. $
 But
$$1-w_c\bar{w}_c=\sqrt{\Delta}\frac{-\sqrt{\Delta}\pm\epsilon_0}{2\epsilon_+\epsilon_-},$$ 
and  in order to assure  $w_c\in\mc{D}_1$, {\it we have to chose the critical  value}
 $ w_c=\frac{-\epsilon_0 +\sqrt{\Delta}}{2\epsilon_+}. $

The solution $\eta_c$  of $\frac{\pa \mc{H}}{\pa \eta}=0$,  equivalent with the
solution of $\dot{x}=0$, $\dot{y}=0$ in (\ref{qqXY}), is  
$\eta_c=
2\frac{2\bar{\epsilon}_a\epsilon_--\epsilon_a\epsilon_0}{\Delta}$. 

We find  the {\it Hessian function}  $H(w,\eta)$ attached to the energy function (\ref{realHH})
around  the critical point $(w_c,\eta_c)$
\begin{equation}\label{hesse}
H(w,\eta)=g  w\bar{w}
+\frac{\epsilon_0}{2}\eta\bar{\eta}+\epsilon_+\eta^2+\epsilon_-\bar{\eta}^2,
{\mbox{\quad where\quad}}
g =\frac{k}{2\sqrt{\Delta}}\left(\frac{4\epsilon_+\epsilon_-}{\epsilon_0-\sqrt{\Delta}}\right)^2 >0
\end{equation}
In general, the critical point  $(w_c,\eta_c)$ of the  energy function
(\ref{realHH}) is non-degenaerate, and the Hessian function
(\ref{hesse}) is positive definite (of index 2)  if $p+2m>0$
(respectively, $p+2m<0$). \\[1ex]

{\bf{Acknowledgement.}} Professor  Pierre Bieliavsky raised me the problem of
  finding the symplectomorphism $FC$  for the K\"ahler two-form on  Siegel-Jacobi domains at the
  XXVIII  Workshop on Geometric Methods in
  Physics in Bia\l owie\.{z}a, Poland. Discussions and correspondence
  with Professor  Pierre Bieliavsky,  Doctor    Yannick
  Voglaire, Professor Rolf Berndt and Professor Jae-Hyun Yang  are warmly acknowledged.
 This investigation  was partially supported by
the CNCSIS-UEFISCSU project PNII- IDEI 454/2009,  Cod
ID-44.

\end{document}